\documentclass{article}

\usepackage{amsmath}
\usepackage{amssymb}
\usepackage{amsthm}

\DeclareMathOperator{\ord}{ord}
\DeclareMathOperator{\length}{length}
\DeclareMathOperator{\gr}{gr}

\DeclareMathOperator{\rad}{rad}
\DeclareMathOperator{\Bl}{Bl}
\DeclareMathOperator{\embdim}{embdim}
\DeclareMathOperator{\Spec}{Spec}
\DeclareMathOperator{\height}{height}

\newcommand{\M}{\ensuremath{\mathfrak{M}}}
\renewcommand{\a}{\ensuremath{\mathfrak{a}}}
\newcommand{\p}{\ensuremath{\mathfrak{p}}}
\newcommand{\V}{\ensuremath{\mathcal{V}}}

\newtheorem{proposition}{Proposition}[section]
\newtheorem{lemma}[proposition]{Lemma}

\theoremstyle{definition}

\newtheorem{example}[proposition]{Example}

\title{On the degree function coefficient of a simple complete ideal in dimension two}
\author{Raymond Debremaeker\\
K.U. Leuven, Celestijnenlaan 200I,\\
3001 Heverlee, Belgium\\
\texttt{raymond.debremaeker@wet.kuleuven.be}
}
\date{}

\begin{document}
\maketitle

\begin{abstract}
Let $(R,\M)$ be a two-dimensional regular local ring with
algebraically closed residue field.  Let $I$ be a simple complete
$\M$-primary ideal of $R$ and let $w$ denote its unique Rees
valuation.  Then the degree function coefficient $d(I,w)=1$.

In this note a short proof of this result is given.

Keywords: Regular local ring, simple complete ideal, degree function
coefficient

2000 Mathematics Subject Classification: 13B22, 13H05, 13H10
\end{abstract}


\section{Introduction}

\noindent

In the theory of complete ideals in two-dimensional regular local
rings the following results holds \cite{debre2}.

\emph{Let $(R,\M)$ be a two-dimensional regular local ring with
algebraically closed residue field $k$.  If $I$ is a simple complete
$\M$-primary of $R$ with unique Rees valuation $w$, then the degree
function coefficient $d(I,w)=1$.}

In the first part of Remark 3.5 in \cite{debre1} this result was
erroneously presented as a consequence of Proposition 3.4 of that
paper.  To correct this we show how the above result can be obtained
from \cite[Propopsition 3.3]{debre1}.  For definitions and
background information the reader is referred to \cite{debre1}.

\begin{proof} First, we recall a few facts from the theory of
complete ideals in two-dimensional regular local rings (with
algebraically closed residue field).\begin{itemize}\item\emph{Every
complete $\M$-primary ideal $I$ of $R$ is normal and minimally
generated.}

Indeed, $I$ is normal since we have in $R$ that the product of
complete ideals is complete again.  Since $R/\M$ is infinite, there
exists an element $x\in\M$ such that $x\V=\M\V$ for all Rees
valuation rings $\V$ of $I$.  Hence $R\Bigl[ \frac{\M}{x} \Bigr]$ is
contained in every Rees valuation ring $\V$ of $I$ and this implies
that \[IR\Bigl[ \frac{\M}{x} \Bigr]\cap R=I,\] i.e., $I$ is
\emph{contracted} from $R\Bigl[ \frac{\M}{x} \Bigr]$.  By
Proposition 2.3 in \cite{Huneke}, we know this is equivalent to
\[\mu(I)=\ord_R(I)+1,\] where $\mu(I)$ denotes the minimal number of
generators of $I$.

In other words
\[\mu(I)=\dim_k\left(\frac{\M^r}{\M^{r+1}}\right)\qquad{\rm with} \:
r:=\ord_R(I),\] i.e., $I$ is \emph{minimally generated} (in the
sense of Definition 3.1 in \cite{debre1}).

\item \emph{Every simple complete $\M$-primary ideal $I$ of $R$ is
quasi-one-fibered.}

If a complete $\M$-primary ideal $I$ of $R$ is simple, then $I$ has
precisely one immediate base point, say $(R_1,\M_1)$.  To see this
we first observe that all immediate base points of $I$ are lying on
the chart $R\Bigl[ \frac{\M}{x} \Bigr]$, where $R\Bigl[ \frac{\M}{x}
\Bigr]$ is contained in every Rees valuation ring of $I$ as in the
previous point.  Next, the transorm $I'$ of $I$ in $R\Bigl[
\frac{\M}{x} \Bigr]$ is simple and complete (see Huneke
\cite[Proposition 3.4 and Proposition 3.5]{Huneke}).  This implies
that $I'$ is contained in just one prime ideal $M_1$ of $R\Bigl[
\frac{\M}{x} \Bigr]$.  Thus $R_1:=R\Bigl[ \frac{\M}{x} \Bigr]_{M_1}$
is the unique immediate base point of $I$ and the transform
$I^{R_1}=I'_{M_1}$ is simple and complete.  Thus every simple
complete $\M$-primary ideal $I$ of $R$ is \emph{quasi-one-fibered}
in the sense of Definition 1.7 in \cite{debre1}.

Hence \[T(I)\subseteq\{v_{\M},w\}\] with $w\in T(I)$ (see
\cite[Proposition 1.5]{debre1}).  Here $T(I)$ denotes the set of
Rees valuations of $I$.

Actually we know that $T(I)=\{w\}$ because of Section 4 in
\cite{Huneke}.

So if $I$ is any simple complete $\M$-primary ideal in a
two-dimensional regular local ring $(R,\M)$ with algebraically
closed residue field, then all the conditions of Proposition 3.3 in
\cite{debre1} are satisfied.\end{itemize} Let
\[(R,\M)<(R_1,\M_1)<\cdots<(R_s,\M_s)\] denote the unique quadratic
sequence determined by the simple complete $\M$-primary ideal $I$ of
$R$ (according to Proposition 1.6 in \cite{debre1}).  By a repeated
use of Proposition 3.3 in \cite{debre1}, it follows that
\[d(I,w)=d(I^{R_1},w)=\ldots=d(I^{R_s},w).\] Since $I^{R_s}=\M_s$
and $w=\ord_{R_s}$-valuation, we have that
\[d(I,w)=d(\M_s,\ord_{R_s}).\] Since $\ord_{R_s}$ is the unique Rees
valuation of $\M_s$, it follows from the theory of degree functions
(see \cite[Theorem 4.3]{ReesSharp}) that
\[e(\M_s)=d(\M_s,\ord_{R_s})\ord_{R_s}(\M_s).\]  As $e(\M_s)=1$, it
follows that $d(\M_s,\ord_{R_s})=1$.  Thus \[d(I,w)=1.\]
\end{proof}


\begin{thebibliography}{00}
\bibitem{debre1} R. Debremaeker, Quasi-one-fibered ideals in two-dimensional Muhly local
domains, J. Algebra 344 (2011) 14--46.
\bibitem{debre2} R. Debremaeker, V. Van Lierde, The effect of quadratic transformations on degree functions,
Beitr\"age Algebra Geom.\ 47 (2006) 121--135.
\bibitem{Huneke} C. Huneke, Complete ideals in two-dimensional regular local rings, in: Commutative Algebra,
Berkeley, CA, 1987, in: Math. Sci. Res. Inst. Publ., vol. 15,
Springer-Verlag New York, 1989, pp. 325--338.
\bibitem{ReesSharp} D. Rees, R.Y. Sharp, On a theorem of B. Teissier on multiplicities of ideals in local rings,
J. London Math.\ Soc.\ 18(2) (1978) 449--463.
\end{thebibliography}
\end{document}